\documentclass[a4paper,10pt,leqno]{amsart}
\usepackage[bulgarian,english]{babel}
\usepackage[cp1251]{inputenc}
\usepackage[T2A]{fontenc}
\usepackage{amsmath}
\usepackage{amssymb}

\usepackage[unicode,pdftex,colorlinks,pagebackref]{hyperref}

\newcommand{\ds}{\displaystyle}
\newcommand{\bp}{\begin{proof}[Proof:]}
\newcommand{\ep}{\end{proof}}

\newtheorem{thm}{Theorem}
\newtheorem{lem}{Lemma}
\newtheorem{Conjecture}{Conjecture}
\newcommand{\R}{\ensuremath{\mathbb{R}}}

\newcommand{\N}{\ensuremath{\mathbb{N}}}
\newcommand{\Z}{\ensuremath{\mathbb{Z}}}

\begin{document}

\title[On the distribution of $\alpha p+\beta$ modulo one]
{On the distribution of $\alpha p+\beta$ modulo one for primes of type $p=ar^2+1$}
\author{T. Todorova}

\date{}

\begin{abstract}
A classical problem in analytic number theory is to study the
distribution of fractional part $\alpha p+\beta$ modulo 1, where $\alpha$ is irrational
and $p$ runs over the set of primes. We consider the subsequence
generated by the primes $p$ such that $p=ar^2+1$
and prove that its distribution has
a similar property.

\smallskip


\end{abstract}

\maketitle

\section{Introduction and statements of the result}
There is an old but hitherto unresolved problem if there
exist infinitely many primes of the form $n^2+ 1$, where $n$ is an integer.
In 1922 G. H. Hardy and J. E. Littlewood gave the
following conjecture:
\begin{Conjecture} Suppose $a,\,b,\,c\in \Z$ with $a>0$,
$GCD(a,\,b,\,c)=1$, $a+b$ and $c$ are not both even, and
$D=b^2-4ac$ is not a square. Let $P_f(x)$ be the number of primes
$p\le x$ of the form $p=f(n)=an^2+bn+c$ with $n\in \Z$.
Then
\begin{equation*}
  P_f(x)\asymp GCD(2,\,a+b)\frac{\mathcal{\sigma}(D){x}}{\sqrt{a}\log x}\prod\limits_{p|a,\,p|b\atop{p>2}}\frac{p}{p-1}\,,
\end{equation*}
where
\begin{equation*}
  \mathcal{\sigma}(D)=\prod\limits_{p\not|a\atop{}p>2}\bigg(1-\frac{\binom{D}{p}}{p-1}\bigg)
\end{equation*}
\end{Conjecture}
One may find
several results on approximations to this problem in the literature.
In 1978 Iwaniec \cite{Iwn2} proved that there are infinitely many numbers of the form $n^2+1$ with at most two prime factors.
In 1952 Ankeny \cite{An} proved that, assuming the extended Riemann hypothesis for $L$-functions on Hecke characters, there are infinitely many primes of the form $x^{2}+y^{2}$ with $y=O(\log x)$.
In 2019 Merikoski \cite{Mer} improving on previous works and showed that there are infinitely many numbers of the form $n^2+1$ with greatest prime factor at least $n^{1.279}$.
The Brun sieve establishes an upper bound on the density of primes having the
form $p=n^{2}+1$ : there are $O({\sqrt {x}}/\log x)$ such primes up to $x$.
It then follows that almost all numbers of the form $n^{2}+1$ are composite.

In 2006, Baier, Zhao \cite{BZ} proved that for $\varepsilon>0$ there exist infinitely many primes of the form $p = am^2+1$ such that
$a\le p^{5/9+\varepsilon}$.
Later Matom\"{a}ki \cite{Mat} improves their result and proved that for $\varepsilon>0$ there exist infinitely many primes of the form $p = aq^2+1$ such that
$a\le p^{1/2+\varepsilon}$ and $q$ is prime.

In the present paper we consider another popular problem with
primes.
Let $\alpha $ be irrational real number, $\beta $ be real and let
$||x||=~\min\limits_{n\in \mathbb{Z}}|x-~n|$. In 1947 Vinogradov \cite{Vin4} proved that if
$0<\theta<1/5$ then there are infinitely many primes
$p$ such that
\begin{equation*}
    ||\alpha p+\beta||<p^{-\theta }.
\end{equation*}
Latter the upper bound for $\theta $ was improved and the
strongest published result is due Matomaki \cite{Mat2} with $\theta<1/3$.

Our result is a hybrid between the two problems mentioned above.
We shall prove the following

\begin{thm} Let $\alpha \in \mathbb{R}\backslash \mathbb{Q}$ ,
$\beta \in \mathbb{R}$
 and let $\,0<\theta< 1/108$. Then for any arbitrary small $\varepsilon>0$ there are
infinitely many primes $p$ such that $p=ar^2+1$ with $a,\,r\in \mathbb{N}$,
\begin{equation*}
  a\le p^{2/3+4\theta+\varepsilon}\quad\hbox{and}\quad ||\alpha p+\beta||<p^{-\theta }\,.
\end{equation*}
\end{thm}
{\bf Acknowledgements:} This work was supported by Sofia University
Scientific Fund, grant 80-10-99/2023.

\section{Notation}

Let $x$ be a sufficiently large real number,
\begin{equation}\label{us2}
    \begin{split}
    \delta &=\delta (x)=x^{-\theta },\quad \,K=\delta ^{-1}\log ^2x\,,\\
    y&=x^{1/6-2\theta-\eta/2 },\quad \theta<\frac{1}{108}\,,
    \end{split}
\end{equation}
where $\eta$ is arbitrary small and positive number, which we choose later.
By $p$ we always denote prime. As usual
$\varphi (n),\,\Lambda(n),\,\tau _k(n)$ are Euler's function, Mangoldt's function and the
number of solutions of the equation $m_1m_2\ldots m_k=n$ in
natural numbers $m_1,\,\ldots,m_k$, $\tau _2(n)=\tau (n)$ and
\begin{equation*}
  \psi(x,\,d,\,a)=\sum\limits_{n\le x\atop{n\equiv a (d)}}\Lambda(n)\,.
\end{equation*}
With $||y||$ we denote
the distance from $y$ to the
nearest integer, $e(y)=e^{2\pi iy}$ and if $X<x\le 2X$ we will write $x\sim X$.
Instead of $m\equiv n\,\pmod {k}$
we write for simplicity $m\equiv n(k)$.
The letter
$\varepsilon$ denotes an arbitrary small positive number, not the
same in all appearances. For example this convention allows us to
write $x^{\varepsilon }\log x\ll x^{\varepsilon }$.

\section{Some Lemmas}
\begin{lem} \label{Tau}
Let $ k,\, l,\,m,\, n \in \mathbb{N} $; $ X, \varepsilon \in \mathbb{R}$;
$X\ge 2$, $k\ge 2$ and $ \varepsilon >0$. Then
\begin{equation*}
\begin{array}{llll}
(\mathrm{i})  \qquad&\ds\sum\limits_{n\le X}  \big( \, \tau_k(n) \, \big)^l  \ll_{ { }_{k, l} }
                         X(\log X)^{k^l-1} \; ;
                         \qquad \qquad  \; \\[12pt]
(\mathrm{ii})  \qquad&\ds\tau_k(n)\ll_{ { }_{k, \varepsilon} } n^{\varepsilon} \; .
\end{array}
\end{equation*}
\end{lem}
\bp See \cite{VinBasic}, ch. 3.
\ep

\begin{lem} \label{Trsum1}
Let $ X\ge 1 $ and $a,\,d\in \mathbb{N}$. Then
\begin{equation*}
   \bigg|
\mathop{\sum}_{\substack{
                      n\le X \\
                      n\equiv a \,(d) }}
e(\alpha n)\bigg|
           \ll \min \left(\frac{X}{d},
                 \frac{1} { ||\alpha d||} \right)
\end{equation*}
\end{lem}
\bp See \cite{Kar}, ch.6, \S 2.
\ep

\begin{lem}\label{BZphi} Let $\varepsilon$ and $A$ be arbitrary positive constants.
If $(\log x)^{A+1}\ll y\ll x^{2/9-\varepsilon}$, then
\begin{equation*}
  \sum\limits_{n\sim x}\Lambda(n+1)\sum\limits_{q\sim y\atop{q^2|n}}1=\frac{x}{2\zeta(2)y}+O\bigg(\frac{x}{y(\log x)^A}\bigg)\,,
\end{equation*}
where the $O$-constant depends only on $\varepsilon$ and $A$.
\end{lem}
\bp See Lemma 10, \cite{BZ}.
\ep


\begin{lem}\label{Mat} Let $x,\,M,\,J\in \R^{+}$, $\mu,\,\zeta\in\N$
and $\alpha\in \mathbb{R}\setminus\mathbb{Q}$ satisfies
conditions
\begin{equation}\label{alfa}
  \bigg|\alpha- \frac{a}{q}\bigg|<\frac{1}{q^2},\quad a \in \mathbb{Z}, \,;\; q \in \mathbb{N }, \quad (a, q) = 1,\quad q \geq 1\,.
\end{equation}
Then for every arbitrary small $\varepsilon>0$ the inequality
\begin{multline*}
  \sum\limits_{m\sim M}\tau_{\mu}(m)\sum\limits_{j\sim J}\tau_{\zeta}(j)
\min\bigg\{\frac{x}{m^2j},\,\frac{1}{||\alpha m^2j||}\bigg\}\\
\ll x^{\varepsilon}\bigg(MJ+\frac{x}{M^{3/2}}+\frac{x}{Mq^{1/2}}+\frac{x^{1/2}q^{1/2}}{M}\bigg)
\end{multline*}
is fulfilled.
\end{lem}
\bp See Lemma 8, \cite{Mato1}.
\ep

\section{Auxiliary results}

\begin{lem}\label{TT} Let $x,\,M,\,J\in \R^{+}$, $\mu,\,\zeta\in \mathbb{N}$
and $\alpha\in \mathbb{R}\setminus\mathbb{Q}$ satisfies
conditions (\ref{alfa}). Then for any $\varepsilon>0$ the inequality
\begin{multline*}
  \sum\limits_{m\sim M}\tau_{\mu}(m)\sum\limits_{j\sim J}\tau_{\zeta}(j)\min\bigg\{\frac{x}{m^4j},\,\frac{1}{||\alpha m^4j||}\bigg\}
\le x^{\varepsilon}\bigg(MJ+\frac{x}{M^{\frac{25}{8}}}+\frac{x}{M^3q^{\frac{1}{8}}}+
\frac{x^{\frac{7}{8}}q^{\frac{1}{8}}}{M^3}\bigg)\,
\end{multline*}
is fulfilled.
\end{lem}
\bp
Our proof is similar to proof of Lemma 8, \cite{Mato1}.
Let
\begin{equation}\label{H}
  H = \dfrac{x}{M^4J}\,.
\end{equation}
If $H \le 2$, then trivially from Lemma \ref{Tau} (ii) we get
\begin{equation}\label{Hsmall}
  G \ll x^{\varepsilon}MJ\,.
\end{equation}
So we can assume that $H > 2$.
From Lemma \ref{Tau} (ii) it is obviously that
\begin{equation*}
  G\ll x^{\varepsilon}\sum\limits_{m\sim M}
\sum\limits_{j\sim J}
  \min\bigg\{\frac{x}{m^4j},\,\frac{1}{||\alpha m^4j||}\bigg\}\,.
\end{equation*}

We apply the Fourier expansion to
function $\min\bigg\{\dfrac{x}{m^4j},\,\dfrac{1}{||\alpha m^4j||}\bigg\}$ and get
\begin{equation*}
  \min\bigg\{\frac{x}{m^4j},\,\frac{1}{||\alpha m^4j||}\bigg\}=
  \sum\limits_{0<|h|\le H^2}w(h)e(\alpha m^4jh) +O(\log x)\,,
\end{equation*}
where
\begin{equation}\label{wh}
  w(h)\ll \min\bigg\{\log H,\,\frac{H}{|h|}\bigg\}\,.
\end{equation}
Then
\begin{equation}\label{GI}
  |G|\ll x^{\varepsilon}\sum\limits_{0<|h|\le H^2}|w(h)|
\sum\limits_{j\sim J}\bigg |\sum\limits_{m\sim M}e(\alpha m^4jh)\bigg |+
  MJ\log x\,.
\end{equation}
So if
\begin{equation*}
  G(H_0)=\sum\limits_{h\sim H_0}\sum\limits_{j\sim J}
  \bigg|\sum\limits_{m\sim M}e(\alpha m^4jh)\bigg|\,.
\end{equation*}
then using (\ref{wh}) we have
\begin{equation}\label{G}
  G\ll x^{\varepsilon}\bigg(MJ+\max\limits_{1\le H_0\le H_1}G(H_0)+
\max\limits_{H_1< H_0\le H^2}\frac{H}{H_0}G(H_0)\bigg)\,.
\end{equation}
We shall evaluate the sum $G(H_0)$. Applying the Cauchy-Schwarz inequality we obtain
\begin{align*}
  G^2(H_0)  &\ll x^{\varepsilon}H_0J \sum\limits_{h\sim H_0}\sum\limits_{j\sim J}
  \bigg|\sum\limits_{m\sim M}e(\alpha m^4jh)\bigg|^2\,\\
&\ll x^{\varepsilon}H_0J \sum\limits_{h\sim H_0}\sum\limits_{j\sim J}
  \sum\limits_{m_1\sim M}\sum\limits_{m_2\sim M}e\big(\alpha (m_2^4-m_1^4)jh\big)\,.
\end{align*}
Substituting $m_2=m_1+t$, where $0\le |t|\le M$ we get
\begin{equation}\label{G2H0}
  G^2(H_0)  \ll x^{\varepsilon}\bigg(H_0^2J^2M+H_0JG_1(H_0)\bigg)\,,
\end{equation}
where
\begin{equation*}
  G_1(H_0)=\sum\limits_{h\sim H_0}\sum\limits_{j\sim J}
  \sum\limits_{0<|t|<M}\bigg|\sum\limits_{m_1\sim M}e\big(\alpha (4m_1^3t+6m_1^2t^2+4m_1t^3)jh\big)\bigg|\,.
\end{equation*}
Applying again the Cauchy-Schwarz inequality we obtain
\begin{multline*}
  G_1^2(H_0) \ll H_0JM
\sum\limits_{h\sim H_0}\sum\limits_{j\sim J}
  \sum\limits_{0<|t|<M}\sum\limits_{m_2\sim M}\\
\times\sum\limits_{m_1\sim M}e\big(\alpha (4(m_2^3-m_1^3)t+6(m_2^2-m_1^2)t^2+4(m_2-m_1)t^3)jh\big)\,.
\end{multline*}
Substituting $m_2=m_1+\ell$, where $0\le |\ell|\le M$ we get
\begin{equation}\label{H0H1}
  G_1^2(H_0) \ll  H_0^2J^2M^3+
H_0JMG_2(H_0)\,,
\end{equation}
where
\begin{equation*}
  G_2(H_0)=\sum\limits_{h\sim H_0}\sum\limits_{j\sim J}
  \sum\limits_{0<|t|<M}\sum\limits_{0<|\ell|<M}\bigg|\sum\limits_{m_1\sim M}e\big(12\alpha (m_1^2t\ell+m_1\ell^2t+m_1\ell t^2)jh\big)\bigg|\,.
\end{equation*}
Applying again the Cauchy-Schwarz inequality we obtain
\begin{multline*}
  G_2^2(H_0) \ll H_0JM^2
\sum\limits_{h\sim H_0}\sum\limits_{j\sim J}
  \sum\limits_{0<|t|<M}\sum\limits_{0<|\ell |<M}\sum\limits_{m_2\sim M}\\
\times\sum\limits_{m_1\sim M}e\big(12\alpha ((m_2^2-m_1^2)t\ell +6(m_2-m_1)\ell ^2t+(m_2-m_1)\ell t^2)jh\big)\,.
\end{multline*}
Substituting $m_2=m_1+z$, where $0\le |z|\le M$ we get
\begin{equation*}
  G_2^2(H_0) \ll  H_0^2J^2M^5+H_0JM^2
\sum\limits_{h\sim H_0}\sum\limits_{j\sim J}
  \sum\limits_{0<|t|<M}\sum\limits_{0<|\ell |<M}\sum\limits_{0<|z|<M}\\
\bigg|\sum\limits_{m_1\sim M}e\big(24\alpha m_1zt\ell jh\big)\bigg |\,.
\end{equation*}
Let $u=24tz\ell zjh$. Then using Lemma \ref{Trsum1} we obtain
\begin{align}\label{GH0square}
  G_2^2(H_0) &\ll H_0^2J^2M^5 +
 H_0JM^2\sum\limits_{u\le 24H_0JM^3}\tau_6(u)\min\bigg\{\frac{H_0JM^4}{u},\,\frac{1}{||\alpha u||}\bigg\}\notag\\
&\ll
x^{\varepsilon}\bigg(H_0^2J^2M^5+
\frac{H_0^2J^2M^6}{q}+H_0JM^2q\bigg)\,.
\end{align}
From (\ref{G2H0}) we obtain
\begin{equation*}
  G(H_0)  \ll x^{\varepsilon}\bigg(H_0JM^{7/8}+
\frac{H_0JM}{q^{\frac{1}{8}}}+
H_0^{\frac{7}{8}}J^{\frac{7}{8}}M^{\frac{1}{2}}q^{\frac{1}{8}}\bigg)\,.
\end{equation*}
Choosing $H_0=H$ from (\ref{G}), (\ref{G2H0}), (\ref{H0H1}),  (\ref{GH0square}),(\ref{H}) and (\ref{Hsmall}) we get
\begin{equation}\label{Gmax}
  G\ll x^{\varepsilon} \bigg(JM+\frac{x}{M^{25/8}}+\frac{x}{M^3q^{1/8}}+
\frac{x^{\frac{7}{8}}q^{\frac{1}{8}}}{M^3}\bigg)\,.
\end{equation}

\ep

\section{Proof of Theorem 1}

Let $s(n)=n/r^2$ , where $n\in \mathbb{N}$ and
$r^2$ is the largest square dividing $n$.
It is clear that $s(n)=1$ if and only if $n$ is a perfect
square. Let $p$ be a prime. Then $s(p-1)\le p^{2/3+4\theta+\eta}$
if and only if there exists $r\in \mathbb{N}$
such that $r^2 |p-1$ and $r\ge (p-1)^{1/6-2\theta-\eta/2}$. Therefore the number $\Gamma (x)$ of primes
$p\sim x-1$ such that $s(p-1)\le p^{2/3+4\theta+\eta
}$ and $||\alpha p+\beta||<p^{-\theta}$ is
\begin{equation*}
  \Gamma (x)\ge \sum\limits_{\substack{ n\sim x-1 \\ ||\alpha n+\beta||<n^{-\theta}}}\frac{\Lambda(n)}{\tau(n-1)}\sum\limits_{r\ge y\atop{r^2|n-1}}1
  -C_2x^{1/2+\varepsilon}\,,
\end{equation*}
where $y=x^{1/6-2\theta-\eta/2}$ and the residual member has been received
from addends of the form $n=s^k,\,k\ge 2$ and $s$-prime.
Using $\tau(n)\ll n^{\varepsilon/2}$ we get
\begin{multline}\label{Gamma}
  \Gamma (x)\ge C_1x^{-\varepsilon/2}\sum\limits_{\substack{ n\sim x-1 \\ ||\alpha n+\beta||<n^{-\theta}}}
  \Lambda(n)\sum\limits_{r\ge y\atop{r^2|n-1}}1
  -C_2x^{1/2+\varepsilon}\\
  \ge C_1x^{-\varepsilon/2}\sum\limits_{\substack{ n\sim x-1 \\ ||\alpha n+\beta||<n^{-\theta}}}\Lambda(n)
  \sum\limits_{r\sim y\atop{r^2|n-1}}1
  -C_2x^{1/2+\varepsilon}\,.
\end{multline}
As in \cite{TT} we take a periodic with period 1 function such that
\begin{equation*}\label{hi}
\begin{aligned}
    0<\chi (t) &<1 \quad \mbox { if }\quad  -\delta< t< \delta;\\
    \chi (t) &=0 \quad \mbox { if }\quad\quad\;  \delta \le t\le 1-\delta,
\end{aligned}
\end{equation*}
and which has a Fourier series
\begin{equation}\label{hi1}
    \chi (t)= \delta +\sum\limits_{|k|>0}c(k)e(kt),
\end{equation}
with coefficients satisfying

\begin{align}\label{hi2}
    c(0)&=\delta,&\nonumber\\
    c(k)&\ll \delta \;\;\;\;\;\mbox{ for all } k,\\
    \sum\limits_{|k|>K}|c(k)|&\ll
    x^{-1}\nonumber
\end{align}
and $\delta $ and $K$ satisfying the conditions (\ref{us2}).
The existence of such a function is a consequence of a well known
lemma of Vinogradov (see \cite{Kar}, ch. 1, \S 2).
Then from (\ref{Gamma}) we get
\begin{equation}\label{GamaGama1}
   \Gamma(x)
    \ge C_1x^{-\varepsilon/2}\Gamma_1(x)-C_2x^{1/2+\varepsilon}\,,
\end{equation}
where
\begin{equation*}
  \Gamma_1(x)=\sum\limits_{n\sim x-1}\Lambda(n)\chi (\alpha n+\beta )\sum\limits_{r\sim y\atop{r^2|n-1}}1\,.
\end{equation*}
From the Fourier expansion (\ref{hi1}) of $\chi(t)$ we get
\begin{equation}\label{RavGama1}
    \Gamma_1 (x)=\delta\bigg(\Gamma_2(x)+\Gamma_3(x)+O\bigg(\frac{1}{y}\bigg)\bigg),
\end{equation}
where
\begin{equation*}
\Gamma_2(x)=\sum\limits_{r\sim y}\sum\limits_{n\sim x \atop{n\equiv 1(r^2)}}\Lambda(n)\quad \hbox{and}\quad
\Gamma_3(x)=\sum\limits_{0<|k|\le K}c(k)
\sum\limits_{r\sim y}\sum\limits_{n\sim x \atop{n\equiv 1(r^2)}}e(\alpha kn)\Lambda(n)\,.
\end{equation*}
Here we have laid $c(k):=c(k)e(\beta k)$.

From Lemma \ref{BZphi} we have
\begin{equation}\label{Gama2}
  \Gamma_2(x)=\frac{x}{2\zeta(2)y}+O\bigg(\frac{x}{y(\log x)^{A}}\bigg)\,.
\end{equation}

\subsection{Estimate of the amount $\Gamma_3(x)$}\label{g2}
$\\$
First we decompose the sum $\Gamma_3(x)$ into $O(\log x)$ sums of
type
\begin{equation*}\label{ccc}
   W=\sum\limits_{k\sim K_0}c(k)
\sum\limits_{r\sim y}\sum\limits_{n\sim x \atop{n\equiv 1(r^2)}}\Lambda(n)e(\alpha kn)\,,
\end{equation*}
where
\begin{equation}\label{OgrK}
  1\le K_0\le K/2\,.
\end{equation}
Then by Vaughan's identity we can decompose the sum $W$ into $O(\log x)$
type I sums
\begin{equation*}
  W_1=\sum\limits_{k\sim K_0}c(k)\sum\limits_{r\sim y}\sum\limits_{m\sim M}a(m)\sum\limits_{\ell\sim L\atop{m\ell\equiv 1\,(r^2)}}e(\alpha m\ell k)\,,
\end{equation*}

\begin{equation*}
  W_1'=\sum\limits_{k\sim K_0}c(k)\sum\limits_{r\sim y}\sum\limits_{m\sim M}a(m)\sum\limits_{\ell\sim L\atop{m\ell\equiv 1\,(r^2)}}e(\alpha m\ell k)\log \ell\,,
\end{equation*}
where $ML\sim x$, $M\le x^{1/3}$ and into $O(\log x)$ type II
sums
\begin{equation*}
  W_2=\sum\limits_{k\sim K_0}c(k)\sum\limits_{r\sim y}\sum\limits_{m\sim M}a(m)\sum\limits_{\ell\sim L\atop{m\ell\equiv 1\,(r^2)}}b(\ell)e(\alpha m\ell k)\,,
\end{equation*}
where $ML\sim x$, $x^{1/3}\le M\le x^{2/3}$ and $a(m)\ll \tau (m)\log m$, $b(l)\ll \tau (l)\log l$.

First we estimate type I sums rather straightforwardly. We have
\begin{equation*}
    W_1\ll x^{\varepsilon}\sum\limits_{k\sim K_0}\sum\limits_{r\sim y}\sum\limits_{m\sim M}\bigg|\sum\limits_{\ell\sim L\atop{m\ell\equiv 1\,(r^2)}}e(\alpha m\ell k)\bigg|
\end{equation*}
As $L>x^{2/3}> y^2\ge r^2$ we get $\ell=f+r^2t$, where $f=f(m,\,r)$. Using Lemma \ref{Trsum1}
we obtain
\begin{align*}
    W_1&\ll x^{\varepsilon}\sum\limits_{k\sim K_0}\sum\limits_{r\sim y}\sum\limits_{m\sim M}\bigg|\sum\limits_{t\sim \frac{L}{r^2}}e(\alpha mtr^2k)\bigg|\\
    &\ll x^{\varepsilon}\sum\limits_{k\sim K_0}\sum\limits_{r\sim y}\sum\limits_{m\sim M}\min\bigg\{\frac{L}{r^2}\,\frac{1}{||\alpha mr^2k||}\bigg\}\\
    &\ll x^{\varepsilon}\sum\limits_{u\sim K_0M}\tau(u)\sum\limits_{r\sim y}\min\bigg\{\frac{xK_0}{ur^2}\,\frac{1}{||\alpha ur^2||}\bigg\}
\end{align*}
Applying Lemma \ref{Mat} and bearing in mind (\ref{OgrK}) we get
\begin{equation}\label{ocW1}
  W_1\ll x^{\varepsilon}\bigg(yx^{1/3}K+\frac{xK}{y^{3/2}}+\frac{xK}{yq^{1/2}}+\frac{x^{1/2}K^{1/2}q^{1/2}}{y}\bigg)
\end{equation}
Now we will estimate the sum $W_2$.
It is enough to considered the case
\begin{equation*}
  x^{1/3}\le L\le x^{1/2}\quad, \quad x^{1/2}\le M\le x^{2/3}\,.
\end{equation*}
We start by
applying the Cauchy-Schwarz inequality, obtaining
\begin{equation}\label{Wsquare}
  W_2^2\le x^{1+\varepsilon}\bigg(xMK^2 +W_{21}\bigg)\,,
\end{equation}
where
\begin{equation*}
  W_{21}=yMK_0\sum\limits_{k\sim K_0}\sum\limits_{r\sim y}\sum\limits_{m\sim M}
  \sum\limits_{\ell_1,\ell_2\sim L\atop{m\ell_1\equiv 1\,(r^2)\atop{m\ell_2\equiv 1\,(r^2)\atop{\ell_\ne\ell_2}}}}
  b(\ell_1)b(\ell_2)e(\alpha m(\ell_1-\ell_2)k)\,.
\end{equation*}
If  $\ell_1\not= \ell_2$ from $m\ell_i\equiv 1\,(r^2)$ follows $\ell_1\equiv \ell_2\,(r^2)$.
So $\ell_1=\ell_2+tr^2$ with $t\le \dfrac{L}{y^2}$ and
using Lemma \ref{Trsum1} we receive
\begin{multline*}
  W_2^2\ll x^{\varepsilon}yMK_0\sum\limits_{k\sim K_0}\sum\limits_{r\sim y}
  \sum\limits_{\ell_2\sim L}\sum\limits_{t\le \frac{L}{y^2}}\bigg|
  \sum\limits_{m\sim M\atop{m\ell_2\equiv 1\,(r^2)}}e(\alpha mtr^2k)
  \bigg|\\
  \ll x^{\varepsilon}xyK_0\sum\limits_{k\sim K_0}\sum\limits_{r\sim y}
  \sum\limits_{t\le \frac{L}{y^2}}
  \min\bigg\{
  \frac{M}{r^2},\,\frac{1}{||\alpha r^4tk||}
  \bigg\}\,.
\end{multline*}
By substitute $u=tk\ll\dfrac{LK_0}{y^2}$ follows
\begin{equation*}
  W_2^2\ll x^{\varepsilon}xyK_0\sum\limits_{u\ll\frac{LK_0}{y^2}}\tau(u)\sum\limits_{r\sim y}
  \min\bigg\{
  \frac{xK_0}{r^4u},\,\frac{1}{||\alpha r^4u||}
  \bigg\}
\end{equation*}
Using Lemma \ref{TT}, (\ref{Wsquare}) and (\ref{OgrK}) we get
\begin{equation}\label{ocW2}
  W_2\ll x^{\varepsilon}\bigg ( x^{\frac{3}{4}}K+
  \frac{xK}{y^{{\frac{17}{16}}}}+\frac{xK}{yq^{{\frac{1}{8}}}}+\frac{x^{{\frac{15}{16}}}K^{{\frac{15}{16}}}q^{{\frac{1}{16}}}}{y}
  \bigg )
\end{equation}
From (\ref{ocW1}), (\ref{ocW2}) and (\ref{g2}) follows
\begin{equation}\label{ocW}
 \Gamma_3\ll x^{\varepsilon}\bigg(x^{\frac{3}{4}}K+
  \frac{xK}{y^{{\frac{17}{16}}}}+\frac{xK}{yq^{{\frac{1}{8}}}}+\frac{x^{{\frac{15}{16}}}K^{{\frac{15}{16}}}q^{{\frac{1}{16}}}}{y}
  \bigg)\,.
\end{equation}
From the denominators $q_1<q_2<...$ of approximation
of irrational $\alpha$ we choose $x_i$ such that the
equality $\dfrac{xK}{y}=q$ is fulfilled. We take $\theta=\frac{1}{108}-2\varepsilon$
and $\eta =8\varepsilon$. So $y=x^{1/6-2\theta-4\varepsilon }$ and get a sequence
\[
   x_1<x_2<x_3 < \dots \to \infty
\]
such that
\begin{equation}\label{ocGama3}
 \Gamma_3\ll \frac{x^{1-\varepsilon}}{y}\,.
\end{equation}
Hence from (\ref{GamaGama1}), (\ref{RavGama1}), (\ref{Gama2}) and (\ref{ocGama3}) we get

\begin{equation*}
  \Gamma\ge C_1x_i^{-\varepsilon/2}\delta\bigg(\frac{x_i}{2\zeta(2)y_i}
  +O\bigg(\frac{x_i^{1-\varepsilon}}{y_i}\bigg)\bigg)
\end{equation*}
and the proof is completed.

\bigskip
\bigskip

\vbox{
\hbox{Faculty of Mathematics and Informatics}
\hbox{Sofia University ``St. Kl. Ohridsky''}
\hbox{5 J.Bourchier, 1164 Sofia, Bulgaria}
\hbox{ }
\hbox{tlt@fmi.uni-sofia.bg}}

\bigskip
\bigskip
\bigskip

\end{document}